 \documentclass[11pt]{article}

 \usepackage{mathrsfs, amssymb}
 \usepackage[dvips]{color}

 \newtheorem{definition}{Definition}[section]
 \newtheorem{hypothesis}{Hypothesis}[section]
 \newtheorem{lemma}{Lemma}[section]
 \newtheorem{proposition}{Proposition}[section]
 \newtheorem{theorem}{Theorem}[section]
 \newtheorem{corollary}{Corollary}[section]

 \def\blemma{\begin{lemma}\sl{}\def\elemma{\end{lemma}}}
 
 \def\btheorem{\begin{theorem}\sl{}\def\etheorem{\end{theorem}}}

 \def\beqlb{\begin{eqnarray}}\def\eeqlb{\end{eqnarray}}
 \def\beqnn{\begin{eqnarray*}}\def\eeqnn{\end{eqnarray*}}

 \def\mbb{\mathbb}
 \def\mbf{\mathbf}
 \def\mcr{\mathscr}

 \def\proof{\noindent\textit{Proof.~~}}\def\qed{\hfill$\Box$\medskip}

 \def\<{\langle}\def\>{\rangle}

\begin{document}

\noindent{Published in: \textit{Journal of Applied Probability}
\textbf{43} (2006), 1: 289--295

\bigskip\bigskip

\noindent{\Large\textbf{A limit theorem of discrete
Galton-Watson}}

\noindent{\Large\textbf{branching processes with immigration}
\footnote{~Supported by NCET and NSFC Grants (No.~10121101 and
No.~10525103).}}

\bigskip

\noindent Zenghu LI \footnote{~Postal address: School of
Mathematical Sciences, Beijing Normal University, Beijing 100875,
P. R. China. E-mail: lizh@bnu.edu.cn} \quad \textit{Beijing Normal
University}

\bigskip

{\narrower{\narrower

\noindent\textbf{Abstract.} We provide a simple set of sufficient
conditions for the weak convergence of discrete Galton-Watson
branching processes with immigration to continuous time and
continuous state branching processes with immigration.

\smallskip

\noindent\textbf{Mathematics Subject Classification (2000):} 60J80

\smallskip

\noindent\textbf{Keywords and phrases:} Markov chain, branching
process with immigration, limit theorem, weak convergence,
generating function.

\smallskip

\noindent\textbf{Abbreviated Title:} Branching processes with
immigration

\par}\par}

\bigskip

\section{Introduction}

\setcounter{equation}{0}

Let $(g,h)$ be a pair of probability generating functions. By a
\textit{discrete time and discrete state Galton-Watson branching
process with immigration} (DBI-process) corresponding to $(g,h)$
we mean a discrete-time Markov chain $\{y(n): n = 0,1,2,\cdots\}$
with state space $\mbb{N} := \{0,1,2,\cdots\}$ and one-step
transition matrix $P(i,j)$ defined by
 \beqlb\label{1.1}
\sum^\infty_{j=0} P(i,j)z^j = g(z)^ih(z), \qquad i=0,1,2,\cdots,~
0\le z\le 1.
 \eeqlb
The intuitive meaning of the process is clear from (\ref{1.1}). In
particular, if $h(z) \equiv 1$, we simply call $\{y(n): n =
0,1,2,\cdots\}$ a \textit{discrete time and discrete state
Galton-Watson branching process} (DB-process).

Kawazu and Watanabe (1971) studied systematically the limit
theorems of DBI-processes. They also characterized completely the
class of the limit processes as \textit{continuous time and
continuous state branching processes with immigration}
(CBI-processes). Let us consider a special class of the
CBI-processes introduced in Kawazu and Watanabe (1971). Suppose
that $R$ is a function on $[0,\infty)$ defined by
 \beqlb\label{1.2}
R(\lambda) = \beta\lambda - \alpha\lambda^2 -
\int_0^\infty\Big(e^{-\lambda u}-1+\frac{\lambda
u}{1+u^2}\Big)\mu(du),
 \eeqlb
where $\beta \in \mbb{R}$ and $\alpha\ge0$ are constants and
$(1\land u^2)\mu(du)$ is a finite measure on $(0,\infty)$, and $F$
is a function on $[0,\infty)$ defined by
 \beqlb\label{1.3}
F(\lambda) = b\lambda + \int_0^\infty(1-e^{-\lambda u})m(du),
 \eeqlb
where $b\ge0$ is a constant and $(1\land u)m(du)$ is a finite
measure on $(0,\infty)$. A Markov process $\{y(t): t\ge0\}$ with
state space $\mbb{R}_+ := [0,\infty)$ is called a
\textit{CBI-process} if it has transition semigroup
$(P_t)_{t\ge0}$ given by
 \beqlb\label{1.4}
\int_0^\infty e^{-\lambda y} P_t(x,dy)
 =
\exp\bigg\{-x\psi_t(\lambda) -
\int_0^tF(\psi_s(\lambda))ds\bigg\}, \qquad \lambda\ge 0,
 \eeqlb
where $\psi_t(\lambda)$ is the unique solution of
 \beqlb\label{1.5}
\frac{d\psi_t}{dt}(\lambda) = R(\psi_t(\lambda)), \qquad
\psi_0(\lambda) = \lambda.
 \eeqlb
Clearly, the transition semigroup $(P_t)_{t\ge0}$ defined by
(\ref{1.4}) is stochastically continuous. In particular, if
$F(\lambda) \equiv 0$, we simply call $\{y(t): t\ge0\}$ a
\textit{continuous time and continuous state branching process}
(CB-process).

A CBI-process is said to be \textit{conservative} if it does not
explode, that is, $\mbf{P}_x\{y(t) < \infty\} = 1$ for every
$t\ge0$ and $x\in \mbb{R}_+$ where $\mbf{P}_x$ denotes the
conditional law given $y(0)=x$. By Kawazu and Watanabe (1971,
Theorem~1.2), the process is conservative if and only if
 \beqnn
\int_{0+} R^*(\lambda)^{-1} d\lambda = \infty,
 \eeqnn
where $R^*(\lambda) = R(\lambda) \vee 0$. (This is a correction to
equation (1.21) of Kawazu and Watanabe (1971).)

Let $\{b_k\}$ and $\{c_k\}$ be sequences of positive numbers such
that $b_k\to \infty$ and $c_k\to \infty$ as $k\to \infty$. Let
$\{y_k(n): n\ge0\}$ be a sequence of DBI-processes corresponding
to the parameters $\{(g_k,h_k)\}$ and assume $y_k(0) = c_k$.
Suppose that for all $t\ge0$ and $\lambda \ge0$ the limits
 \beqlb\label{1.6}
\lim_{k\to\infty} g_k^{[kt]}(e^{-\lambda/b_k})^{c_k} =
\phi_1(t,\lambda)
 \quad\mbox{and}\quad
\lim_{k\to\infty} \prod_{j=0}^{[kt]-1}h_k(g_k^j(e^{-\lambda/b_k}))
= \phi_2(t,\lambda)
 \eeqlb
exist and the convergence is locally uniform in $\lambda\ge0$ for
each fixed $t\ge0$, where $g_k^j$ denotes the $j$-order
composition of $g_k$ and $[kt]$ denotes the integer part of $kt$.
The following result was proved in Kawazu and Watanabe (1971,
Theorem~2.1):

\btheorem\label{t1.1} Suppose that (\ref{1.6}) holds and
$\phi_1(t,\lambda) <1$ for some $t>0$ and $\lambda >0$. Then
$\{y_k([kt])/b_k: t\ge0\}$ converges in finite-dimensional
distributions to a stochastically continuous and conservative
CBI-process $\{y(t): t\ge0\}$ with transition semigroup given by
(\ref{1.4}). \etheorem

Based on this theorem, Kawazu and Watanabe (1971) showed that,
given each stochastically continuous and conservative CBI-process
$\{y(t): t\ge0\}$, there is a sequence of positive numbers
$\{b_k\}$ with $b_k\to \infty$ and a sequence of DBI-processes
$\{y_k(n): n\ge0\}$ such that $\{y_k([kt])/b_k: t\ge0\}$ converges
in finite-dimensional distributions to $\{y(t): t\ge0\}$. Their
results have become the basis of many studies of branching
processes with immigration; see e.g.\ Pitman and Yor (1982) and
Shiga and Watanabe (1973). On the other hand, since condition
(\ref{1.6}) involves complicated compositions of the probability
generating functions, it is some times not so easy to verify. In
view of the characterizations (\ref{1.1}), (\ref{1.4}) and
(\ref{1.5}) of the two classes of processes, one naturally expect
some simple sufficient conditions for the convergence of the
DBI-processes to the CBI-processes given in terms of the
parameters $(g,h)$ and $(R,F)$. The purpose of this note is to
provide a set of conditions of this type. For the convenience of
proof, we shall discuss the convergence of $\{y_k([\gamma_kt])/k:
t\ge0\}$ for some sequence of positive numbers $\{\gamma_k\}$ with
$\gamma_k\to \infty$, which is slightly different from the scaling
of Kawazu and Watanabe (1971). Instead of the convergence of
finite-dimensional distributions, we shall consider the weak
convergence on the space of c\`adl\`ag functions $D([0,\infty),
\mbb{R}_+)$.

\section{The limit theorem}

\setcounter{equation}{0}

In this section, we prove a limit theorem for DBI-processes on the
space $D([0,\infty), \mbb{R}_+)$. Let $F$ be defined by
(\ref{1.3}). For simplicity we assume the function $R$ is given by
 \beqlb\label{2.1}
R(\lambda) = \beta\lambda - \alpha\lambda^2 -
\int_0^\infty(e^{-\lambda u}-1+\lambda u)\mu(du), \qquad
\lambda\ge 0,
 \eeqlb
where $\beta \in \mbb{R}$ and $\alpha\ge0$ are constants and
$(u\land u^2)\mu(du)$ is a finite measure on $(0,\infty)$. Suppose
that $\{y(t): t\ge0\}$ is a CBI-process corresponding to $(R,F)$.
Let $\{y_k(n): n\ge0\}$ be a sequence of DBI-processes
corresponding to the parameters $\{(g_k,h_k)\}$ and let
$\{\gamma_k\}$ be a sequence of positive numbers. For $0\le
\lambda\le k$ set
 \beqlb\label{2.2}
F_k(\lambda) = \gamma_k[1 - h_k(1-\lambda/k)]
 \eeqlb
and
 \beqlb\label{2.3}
R_k(\lambda) = k\gamma_k[(1-\lambda/k) - g_k(1-\lambda/k)].
 \eeqlb
Let us consider the following set of conditions:

(2.A) As $k\to \infty$, we have $\gamma_k \to\infty$ and
$\gamma_k/k \to$ some $\gamma_0\ge 0$.

(2.B) As $k\to \infty$, the sequence $\{F_k\}$ defined by
(\ref{2.2}) converges to a continuous function.

(2.C) The sequence $\{R_k\}$ defined by (\ref{2.3}) is uniformly
Lipschitz on each bounded interval and converges to a continuous
function as $k\to \infty$.

We remark that conditions (2.B) and (2.C) parallel the sufficient
conditions for the convergence of continuous-time and discrete
state branching processes with immigration, see e.g., Li (1992)
for the discussions in the setting of measure-valued processes.
Based the results of Li (1991), the following lemma can be proved
by modifying the arguments of the proofs of Li (1992, Lemmas~3.4
and~4.1).

\blemma\label{l2.1} (i) Under conditions (2.B) and (2.C), the
limit functions $F$ and $R$ of $\{F_k\}$ and $\{R_k\}$ have
representations (\ref{1.3}) and (\ref{2.1}), respectively. (ii)
For any $(F,R)$ given by (\ref{1.3}) and (\ref{2.1}), there are
sequences $\{\gamma_k\}$ and $\{(g_k, h_k)\}$ as above such that
(2.A), (2.B) and (2.C) hold with $F_k \to F$ and $R_k \to R$.
\elemma

For $\lambda\ge0$ we set
 \beqlb\label{2.4}
S_k(\lambda) = k\gamma_k[(1-\lambda/k) - g_k(e^{-\lambda/k})].
 \eeqlb

\blemma\label{l2.2} Under conditions (2.A) and (2.C), let $R =
\lim_{k\to\infty} R_k$. Then we have
 \beqlb\label{2.5}
\lim_{k\to\infty} S_k(\lambda) = R(\lambda) - \gamma_0\lambda^2/2
 \quad\mbox{and}\quad
\lim_{k\to\infty} \gamma_k[1 - g_k(e^{-\lambda/k})] =
\gamma_0\lambda
 \eeqlb
uniformly on each bounded interval. \elemma

\proof By mean-value theorem we have
 \beqlb\label{2.6}
S_k(\lambda)
 =
R_k(\lambda) - k\gamma_kg_k^\prime(\eta_k)
(e^{-\lambda/k}-1+\lambda/k),
 \eeqlb
where $1-\lambda/k < \eta_k < e^{-\lambda/k}$ and $g_k^\prime$
denotes the derivative of $g_k$. Under condition (2.C), the
sequence $R_k^\prime(\lambda) = \gamma_k[g_k^\prime
(1-\lambda/k)-1]$ is uniformly bounded on each bounded interval
$\lambda \in [0,l]$ for $l \ge 0$. Then $g_k^\prime(1-\lambda/k)
\to 1$ uniformly on each bounded interval. In particular, we have
$g_k^\prime (\eta_k) \to 1$ and the first equality in (\ref{2.5})
follows from (2.A) and (\ref{2.6}). The second equality follows by
a similar argument. \qed

\btheorem\label{t2.1} Suppose conditions (2.A), (2.B) and (2.C)
hold with $F = \lim_{k\to\infty} F_k$ and $R = \lim_{k\to\infty}
R_k$. If $y_k(0)/k$ converges in distribution to $y(0)$, then
$\{y_k([\gamma_kt])/k: t\ge0\}$ converges in distribution on
$D([0,\infty), \mbb{R}_+)$ to the CBI-process $\{y(t): t\ge0\}$
corresponding to $(R,F)$ with initial value $y(0)$. \etheorem

\proof Let $(P_t)_{t\ge0}$ denote the transition semigroup of the
CBI-process corresponding to $(R,F)$. For $\lambda>0$ and $x\ge0$
set $e_\lambda(x) = e^{-\lambda x}$. We denote by $D_1$ the linear
hull of $\{e_\lambda: \lambda>0\}$. Then $D_1$ is an algebra which
strongly separates the points of $\mbb{R}_+$. Let $C_0(\mbb{R}_+)$
be the space of continuous functions on $\mbb{R}_+$ vanishing at
infinity. By the Stone-Weierstrass theorem, $D_1$ is dense in
$C_0(\mbb{R}_+)$ for the supremum norm; see, e.g., Hewitt and
Stromberg (1965, pp.98-99). For $\lambda>0$ we set
 \beqlb\label{2.7}
Ae_\lambda(x)
 =
-e^{-\lambda x}\left[xR(\lambda) + F(\lambda)\right], \qquad x\in
\mbb{R}_+,
 \eeqlb
and extend the definition of $A$ to $D_1$ by linearity. Then $A :=
\{(f,Af): f \in D_1\}$ is a linear subspace of $C_0(\mbb{R}_+)
\times C_0(\mbb{R}_+)$. Since $D_1$ is invariant under
$(P_t)_{t\ge0}$, it is a core of $A$; see, e.g., Ethier and Kurtz
(1986, p.17). With those observations it is not hard to see that
the semigroup $(P_t)_{t\ge0}$ is generated by the closure of $A$;
see, e.g., Ethier and Kurtz (1986, p.15 and p.17). Note that
$\{y_k (n)/k: n\ge0\}$ is a Markov chain with state space $E_k :=
\{0,1/k,2/k,\cdots\}$ and one-step transition probability
$Q_k(x,dy)$ determined by
 \beqnn
\int_{E_k}e^{-\lambda y}Q_k(x,dy)
 =
g_k(e^{-\lambda/k})^{kx}h_k(e^{-\lambda/k}).
 \eeqnn
Then the (discrete) generator $A_k$ of $\{y_k ([\gamma_kt])/k:
t\ge0\}$ is given by
 \beqnn
A_ke_\lambda(x)
 &=&
\gamma_k\Big[g_k(e^{-\lambda/k})^{kx}h_k(e^{-\lambda/k})
- e^{-\lambda x}\Big] \\
 &=&
\gamma_k\Big[\exp\{xk\alpha_k(\lambda)(g_k(e^{-\lambda/k})-1)\}
\exp\{\beta_k(\lambda)(h_k(e^{-\lambda/k})-1)\} - e^{-\lambda
x}\Big],
 \eeqnn
where
 \beqnn
\alpha_k(\lambda)
 =
(g_k(e^{-\lambda/k})-1)^{-1} \log g_k(e^{-\lambda/k})
 \eeqnn
and $\beta_k(\lambda)$ is defined by the same formula with $g_k$
replaced by $h_k$. Under conditions (2.A), (2.B) and (2.C), it is
easy to show that
 \beqnn
\lim_{k\to\infty}(g_k(e^{-\lambda/k})-1) =
\lim_{k\to\infty}(h_k(e^{-\lambda/k})-1) = 0
 \eeqnn
and
 \beqnn
\lim_{k\to\infty} \alpha_k(\lambda) = \lim_{k\to\infty}
\beta_k(\lambda) =1.
 \eeqnn
Then we have
 \beqlb\label{2.8}
A_ke_\lambda(x)
 =
-e^{-\lambda x}\left[x\alpha_k(\lambda)S_k(\lambda) +
x\gamma_k(\alpha_k(\lambda)-1)\lambda + H_k(\lambda)\right] +
o(1),
 \eeqlb
where
 \beqnn
H_k(\lambda)
 =
\gamma_k\beta_k(\lambda)(1-h_k(e^{-\lambda/k})).
 \eeqnn
By elementary calculations we find that
 \beqnn
\alpha_k(\lambda)
 =
1 + \frac{1}{2}(1-g_k(e^{-\lambda/k})) + o(1-g_k(e^{-\lambda/k})),
 \eeqnn
and so $\lim_{k\to\infty} \gamma_k(\alpha_k(\lambda)-1) =
\gamma_0\lambda/2$ by Lemma~\ref{l2.2}. It follows that
 \beqnn
\lim_{k\to\infty}\left[\alpha_k(\lambda)S_k(\lambda) +
\gamma_k(\alpha_k(\lambda)-1)\lambda\right] = R(\lambda).
 \eeqnn
By the argument of the proof of Lemma~\ref{l2.2} one can show that
 \beqnn
\lim_{k\to\infty}H_k(\lambda) = \lim_{k\to\infty} F_k(\lambda) =
F(\lambda).
 \eeqnn
In view of (\ref{2.7}) and (\ref{2.8}) we get
 \beqnn
\lim_{k\to\infty} \sup_{x\in E_k}\left|A_ke_\lambda(x) -
Ae_\lambda(x)\right| = 0
 \eeqnn
for each $\lambda>0$. This clearly implies that
 \beqnn
\lim_{k\to\infty} \sup_{x\in E_k}\left|A_kf(x) - Af(x)\right| = 0
 \eeqnn
for each $f\in D_1$. By Ethier and Kurtz (1986, p.226 and
pp.233-234) we find that $\{y_k([\gamma_kt])/k: t\ge0\}$ converges
in distribution on $D([0,\infty), \mbb{R}_+)$ to the CBI-process
corresponding to $(R,F)$. \qed

By Lemma~\ref{l2.1} and Theorem~\ref{t2.1}, for any functions
$(R,F)$ given by (\ref{1.3}) and (\ref{2.1}), there is a sequence
of positive numbers $\{\gamma_k\}$ and a sequence of DBI-processes
$\{y_k(n): n\ge0\}$ such that $\{y_k([\gamma_kt])/k: t\ge0\}$
converges in distribution on $D([0,\infty), \mbb{R}_+)$ to the
CBI-process corresponding to $(R,F)$.

\section{Generalized Ray-Knight theorems}

\setcounter{equation}{0}

As an example of the applications of their limit theorems, Kawazu
and Watanabe (1971) reproved the Ray-Knight theorems of diffusion
characterizations of the Brownian local time. In this section, we
generalize the results to the case of a Brownian motion with
drift. We refer the reader to Le Gall and Le Jan (1998) for
another adequate formulation of the Ray-Knight theorems for
general L\'evy processes.

Let $A = \alpha d^2/dx^2 + \beta d/dx$ for given constants $\alpha
>0$ and $\beta \in \mbb{R}$. Then $A$ generates a one-dimensional
Brownian motion with drift $(X_t, \mcr{F}_t, \mbf{P}_x)$. The
\textit{local time} of $\{X_t: t\ge0\}$ is a continuous
two-parameter process $\{l(t,x): t\ge0, x\in \mbb{R}\}$ such that
the following property holds almost surely:
 \beqlb\label{3.1}
2\int_B l(s,x)dx = \int_0^t 1_B(X_s)ds, \qquad B\in
\mcr{B}(\mbb{R}),
 \eeqlb
where $\mcr{B}(\mbb{R})$ denote the Borel $\sigma$-algebra of
$\mbb{R}$ and $1_B$ denotes the indicator function of $B$. For
fixed $a\ge0$ let
 \beqlb\label{3.2}
l^{-1}(u,a) = \inf\{t\ge0: l(t,a)=u\}.
 \eeqlb

\btheorem\label{t3.1} The process
 \beqlb\label{3.3}
\xi_u(t) = l(l^{-1}(u,a),a+t), \qquad t\ge0,
 \eeqlb
is a diffusion generated by
 \beqlb\label{3.4}
x\frac{d^2}{dx^2} + \frac{\beta}{\alpha} x\frac{d}{dx}.
 \eeqlb
\etheorem

\proof We follow the ideas of Kawazu and Watanabe (1971,
Example~2.2). For $c\in \mbb{R}$ let $\sigma_c = \inf\{t\ge0:
X_t=c\}$. Let $\delta>0$ and let $u_\delta(x) =
\mbf{P}_x\{\sigma_\delta < \sigma_{-\delta}\} = 1 -
\mbf{P}_x\{\sigma_\delta > \sigma_{-\delta}\}$. Then
$u_\delta(\cdot)$ satisfies
 \beqnn
\alpha \frac{d^2}{dx^2}u_\delta(x) + \beta \frac{d}{dx}u_\delta(x)
= 0, \qquad |x| \le \delta,
 \eeqnn
with $u_\delta(\delta) = 1$ and $u_\delta(-\delta) = 0$. Solving
this boundary value problem we find that
 \beqnn
u_\delta(x) = \frac{\exp\{\beta\delta/\alpha\} - \exp\{-\beta
x/\alpha\}} {\exp\{\beta\delta/\alpha\} -
\exp\{-\beta\delta/\alpha\}}.
 \eeqnn
By a $\delta$-downcrossing at $x\in \mbb{R}$ before time $T>0$ we
mean an interval $[u,v] \subset [0,T)$ such that $X_u = x+\delta$,
$X_v = x$ and $x<X_t<x+\delta$ for all $u<t<v$. Let $\eta_\delta$
denote the number of $\delta$-downcrossings at $0$ before time
$\sigma_{-\delta}$. By the property of independent increments of
the Brownian motion with drift we have
 \beqnn
\mbf{E}_0[z^{\eta_\delta}] = \sum_{i=0}^\infty (1-p)(pz)^i =
\frac{q}{1-pz},
 \eeqnn
where $p=u_{\delta}(0)$, $q=1-p$ and $\mbf{E}_0$ denotes the
expectation under $\mbf{P}_0$. Let $x_i = a + i/k$ for $i\ge0$ and
$k\ge1$ and let $Z_k(i)$ denote the number of $1/k$-downcrossings at
$x_i$ before time $l^{-1}(u,a)$. It is easy to see that $Z_k(i+1)$
is the sum of $Z_k(i)$ independent copies of $\eta_{1/k}$. Thus
$\{Z_k(i): i=0,1,\cdots\}$ is a DB-process corresponding to the
generating function
 \beqnn
g_k(z) = \frac{q_k}{1-p_kz},
 \eeqnn
where $p_k=u_{1/k}(0)$ and $q_k=1-p_k$. By a standard result for
local times of diffusion processes,
 \beqnn
\lim_{k\to 0} Z_{1/k}([kt])/k = l(l^{-1}(u,a),a+t) = \xi_u(t);
 \eeqnn
see It\^o and McKean (1965, p.48 and p.222). Then
Theorem~\ref{t2.1} implies that the limit $\{\xi_u(t): t\ge0\}$ is
a CB-process corresponding to
 \beqnn
R(\lambda) = \lim_{k\to\infty} k^2[(1-\lambda/k) -
g_k(1-\lambda/k)] = \frac{\beta}{\alpha}\lambda - \lambda^2.
 \eeqnn
This proves the desired result. \qed

Kawazu and Watanabe (1971, Theorem~2.3 and Example~2.2) proved the
results of Theorem~\ref{t3.1} in the special case $\beta=0$. In
that case the generating function $g_k$ is actually independent of
$k\ge 1$. In the general case, it seems difficult to check
condition (\ref{1.6}) for the sequence $\{g_k\}$. By similar
arguments as the above we obtain the following

\btheorem\label{t3.2} The process
 \beqlb\label{3.5}
\eta_u(t) = l(l^{-1}(u,a),a-t), \qquad 0\le t\le a.
 \eeqlb
is a diffusion generated by
 \beqlb\label{3.6}
x\frac{d^2}{dx^2} + \frac{\beta}{\alpha} x\frac{d}{dx} +
\frac{d}{dx}.
 \eeqlb
\etheorem

\bigskip

\textbf{Acknowledgement.} I would like to thank the referee for a
number of helpful comments.

\end{document}